\documentclass[10pt,twocolumn,letterpaper,twoside]{article}

\usepackage[letterpaper,left=1in,right=1in,top=1in,bottom=1in]{geometry}

\title{On the $\K$-theory of $\bZ/p^n$ -- announcement}
\author{Benjamin Antieau\thanks{Northwestern University}, Achim
Krause\thanks{Universit\"at M\"unster}, and Thomas
Nikolaus\footnotemark[2]}
\date{April 7, 2022}

\usepackage[pdfstartview=FitH,
            pdfauthor={},
            pdftitle={2021 notes},
            colorlinks,
            linkcolor=reference,
            citecolor=citation,
            urlcolor=e-mail,
            ]{hyperref}

\usepackage{multirow}
\usepackage{amsmath}
\usepackage{amscd}
\usepackage{amsbsy}
\usepackage{amssymb}
\usepackage{verbatim}
\usepackage{eufrak}
\usepackage{eucal}
\usepackage{microtype}
\usepackage{hyperref}
\usepackage{mathrsfs}
\usepackage{amsthm}
\usepackage{stmaryrd}
\usepackage{xy}
\input xy
\xyoption{all}
\usepackage{tikz}
\usetikzlibrary{matrix,arrows}
\usepackage{tikz-cd}
\usepackage{dsfont}
\usepackage[T1]{fontenc}

\usepackage{fancyhdr}
\newcommand{\stackspace}{2.5}
\newcommand{\stack}[2][1cm]{\;\tikz[baseline, yshift=.65ex]%
    {\foreach \k [evaluate=\k as \r using (.5*#2+.5-\k)*\stackspace] in {1,...,#2}{%
    \ifodd\k{\draw[->](0,\r pt)--(#1,\r pt);}%
    \else{\draw[<-](0,\r pt)--(#1,\r pt);}\fi
    }}\;}

\usepackage{makeidx}
\makeindex

\usepackage[bbgreekl]{mathbbol}
\usepackage{amsfonts}
\DeclareSymbolFontAlphabet{\mathbb}{AMSb} 
\DeclareSymbolFontAlphabet{\mathbbl}{bbold}
\newcommand{\Prism}{{ \mathbbl{\Delta}}}

\usepackage{color}
\definecolor{todo}{rgb}{1,0,0}
\definecolor{conditional}{rgb}{0,1,0}
\definecolor{e-mail}{rgb}{0,.40,.80}
\definecolor{reference}{rgb}{.20,.60,.22}
\definecolor{mrnumber}{rgb}{.80,.40,0}
\definecolor{citation}{rgb}{0,.40,.80}

\let\oldmarginpar\marginpar
\renewcommand\marginpar[1]{\-\oldmarginpar[\raggedleft\footnotesize #1]%
{\raggedright\footnotesize #1}}

\newcommand{\Fscr}{\mathcal{F}}

\newcommand{\Nscr}{\mathcal{N}}
\newcommand{\Oscr}{\mathcal{O}}

\newcommand{\E}{\mathrm{E}}
\newcommand{\F}{\mathrm{F}}

\renewcommand{\H}{\mathrm{H}}

\newcommand{\K}{\mathrm{K}}

\newcommand{\R}{\mathrm{R}}

\newcommand{\bF}{\mathbf{F}}

\newcommand{\bQ}{\mathbf{Q}}

\newcommand{\bZ}{\mathbf{Z}}

\newcommand{\syn}{\mathrm{syn}}

\newcommand{\can}{\mathrm{can}}

\renewcommand{\geq}{\geqslant}
\renewcommand{\leq}{\leqslant}

\newcommand{\TP}{\mathrm{TP}}
\newcommand{\TC}{\mathrm{TC}}

\newcommand{\TR}{\mathrm{TR}}

\newcommand{\Prismhat}{\widehat{\Prism}}

\newcommand{\et}{\mathrm{\acute{e}t}}

\DeclareMathOperator{\Spec}{Spec}
\DeclareMathOperator{\Spf}{Spf}

\newcommand{\we}{\simeq}
\newcommand{\iso}{\cong}

\theoremstyle{plain}
\newtheorem{theorem}{Theorem}[section]
\newtheorem*{theorem*}{Theorem}

\newtheorem{proposition}[theorem]{Proposition}

\newtheorem{corollary}[theorem]{Corollary}
\newtheorem*{corollary*}{Corollary}

\theoremstyle{plain}

\theoremstyle{definition}

\newtheoremstyle{named}{}{}{\itshape}{}{\bfseries}{.}{.5em}{#1 \thmnote{#3}}
\theoremstyle{named}

\theoremstyle{definition}

\newtheorem{example}[theorem]{Example}
\newtheorem*{example*}{Example}

\newtheorem*{question*}{Question}

\begin{document}

\maketitle

Quillen introduced higher algebraic $\K$-theory in~\cite{quillen-higher} and computed
the $\K$-groups $\K_*(\bF_q)$ in~\cite{Qui}.  Except in low degrees, the
computation of the $\K$-groups of closely related rings, for example $\bZ/4$, has remained out of
reach. In this paper, we announce new methods for computations of $\K$-groups
of such rings and outline new results. A full account will be given
in~\cite{akn}.

We are interested in rings of the form $\Oscr_K/\varpi^n$ where
$K$ is a finite extension of $\bQ_p$ of degree $d$, $\Oscr_K$ is its ring of integers, and
$\varpi^n$ is the $n$th power of a uniformizer $\varpi$. In particular, $p\in(\varpi^e)$ where $e$ is the
degree of ramification of $K$ over $\bQ_p$. When $n=1$, $\Oscr_K/\varpi^n$ is
the residue field $k=\bF_q$ of $\Oscr_K$, where $q=p^f$ for some $f$,
called the residual degree of the extension.

The problem of computing the $\K$-groups of such rings, and of finite rings in
general, was raised by Swan in the Battelle proceedings~\cite[Prob.~20]{gersten-problems}.

\section{History}

For any field $k$, $\K_0(k)\iso\bZ$ and $\K_1(k)\iso k^\times$.
Quillen showed in~\cite{Qui} that if $\bF_q$ is the finite field with $q=p^f$ elements, then
for $r\geq 1$,
$$\K_r(\bF_q)\iso\begin{cases}
    0&\text{if $r$ is even and}\\
    \bZ/(q^i-1)&\text{if $r=2i-1$}.
\end{cases}$$
Note in particular that there is no $p$-torsion in the $\K$-groups of $\bF_q$.

For each prime $\ell$ and ring $R$, $\K(R,\bZ_\ell)$ denotes the
$\ell$-completion of the $\K$-theory spectrum of $R$. In the main case of interest
to us, namely when $R=\Oscr_K/\varpi^n$, $\K_r(R)$ is finitely generated
torsion for $r>0$ and $\K_r(R,\bZ_\ell)$ is the subgroup of $\ell$-primary
torsion in $\K_r(R)$. 

Gabber's rigidity theorem~\cite{gabber} implies that if
$R$ is a commutative ring which is henselian with respect to an ideal $I$ and
if $\ell$ is invertible in $R$, then $$\K(R;\bZ_\ell)\we\K(R/I;\bZ_\ell).$$
Examples of such henselian pairs are the rings of integers $\Oscr_K$ as above with the ideal
$(\varpi)$ or the quotients $\Oscr/\varpi^n$, again with the ideal
$(\varpi)$. It follows that for $\ell\neq p$ we have
$$\K_*(\Oscr;\bZ_\ell)\iso\K_*(\Oscr/\varpi^n;\bZ_\ell)\iso\K_*(\bF_q;\bZ_\ell)$$
so that these $\ell$-adic $\K$-groups are all determined by Quillen's computation.

The situation of the $p$-adic $\K$-theory of $\Oscr_K$ or $\Oscr_K/\varpi^n$ is very different.
A result of Dundas--Goodwillie--McCarthy~\cite{dundas-goodwillie-mccarthy}
implies that $\K(\Oscr/\varpi^n;\bZ_p)\we \tau_{\geq 0}\TC(\Oscr/\varpi^n;\bZ_p)$,
while work of Hesselholt--Madsen~\cite{hesselholt-madsen-finite} and of
Panin~\cite{panin} implies that $\K(\Oscr_K;\bZ_p)\we\tau_{\geq 0}\TC(\Oscr_K;\bZ_p)$.
Here, $\TC(\Oscr_K;\bZ_p)$ and $\TC(\Oscr_K/\varpi^n;\bZ_p)$
denote the $p$-adic topological cyclic homology spectra of
$\Oscr_K$ and $\Oscr_K/\varpi^n$, respectively. This theory is built from topological
Hochschild homology and is closely connected to $p$-adic cohomology theories
thanks to the work of~\cite{bms2}.
These results make the $p$-adic $\K$-groups amenable to calculation using
so-called trace methods.

Hesselholt and Madsen determine the structure of $\TC_*(\Oscr_K;\bZ_p)\iso\K_*(\Oscr_K;\bZ_p)$ in~\cite{hesselholt-madsen} and
thereby verify the Quillen--Lichtenbaum conjecture for $\Oscr_K$. This
conjecture now follows in general from the proof of the Bloch--Kato conjecture
due to Rost and Voevodsky; see for example~\cite{haesemeyer-weibel}, although
the $p$-adic ranks of the groups $\K_*(\Oscr_K;\bZ_p)$ had previously
been computed by Wagoner~\cite{wagoner}.

The Hesselholt--Madsen approach uses logarithmic de Rham--Witt forms and
$\TR$, i.e., the classical approach to trace method computations. These have
recently been revisited by Liu--Wang~\cite{liu-wang} who describe $\K_*(\Oscr_\K;\bF_p)$, the
$\K$-groups with mod $p$ coefficients, using new cyclotomic
techniques from~\cite{bms2,nikolaus-scholze}.

The result is that $$\K_r(\Oscr_K;\bZ_p)\iso\begin{cases}
    \bZ_p&\text{if $r=0$,}\\
    \H^1_\et(\Spec K,\bZ_p(i))&\text{if $r=2i-1$, and}\\
    \H^2_\et(\Spec K,\bZ_p(i))&\text{if $r=2i-2$,}
\end{cases}$$
where $\bZ_p(i)$ is the $i$th Tate twist.
These cohomology groups are determined by Iwasawa theory: for $i>0$,
\begin{align*}
    \H^1_\et(\Spec K,\bZ_p(i))&\iso\bZ_p^d\oplus\bZ/w_i,\\
    \H^2_\et(\Spec K,\bZ_p(i))&\iso\bZ/w_{i-1},
\end{align*}
where $d$ is the degree of $K$ over $\bQ_p$ and where
$w_i$ is the largest $p$th power $p^\nu$ such that the exponent of the
cyclotomic Galois group
$\mathrm{Gal}(K(\mu_{p^\nu})/K)$ divides $i$. The number $w_i$ is the $p$-part
of a number introduced by Harris--Segal~\cite{harris-segal}, Quillen, and
Lichtenbaum in the setting of the Quillen--Lichtenbaum conjecture. See Weibel's
book~\cite[Chap.~VI]{weibel-kbook} for more details.

Much less is known about the $\K$-theory of the intermediate rings
$\Oscr_K/\varpi^n$ for $1<n<\infty$.
As for fields, $\K_0(\Oscr_K/\varpi^n)\iso\bZ$ and $\K_1(\Oscr_K/\varpi^n)$ is isomorphic to the group of units in
$\Oscr_K/\varpi^n$.

In~\cite{dennis-stein}, Dennis and Stein determined the structure of $\K_2(\Oscr_K/\varpi^n)$.
No other work we are aware of has addressed the $\K$-groups of general rings of the
form $\Oscr_K/\varpi^n$.

In special situations, more is known. First, every ring
$\bF_q[z]/z^n$ is of the form $\Oscr_K/\varpi^n$ for $K$ of ramification degree
at least $n$. The algebraic $\K$-groups
of these truncated polynomial rings have been studied by
Hesselholt--Madsen in~\cite{hesselholt-madsen-truncated} using classical trace
method techniques, by Speirs in~\cite{speirs} using the new approach to
$\TC$ due to Nikolaus--Scholze~\cite{nikolaus-scholze}, and by Sulyma
in~\cite{sulyma} using the approach to $\TC$ via syntomic cohomology due to
Bhatt--Morrow--Scholze~\cite{bms2} and as outlined by Mathew in~\cite{mathew-survey}.

Second, for unramified extension there are some results in low degrees. In the
unramified case, where $e=1$, $\Oscr_K$ is the ring $W(\bF_q)$ of $p$-typical Witt vectors
of the residue field. Brun~\cite{brun} determined the $\K$-groups of $\bZ/p^n$
(i.e., when $e=1$ and $f=1$) up to degree $p-3$ and
Angeltveit~\cite{angeltveit} determined the $\K$-groups of
$W_n(\bF_q)=W(\bF_q)/\varpi^n=W(\bF_q)/p^n$ up
to degree $2p-2$.

Angeltveit also proved an important quantitative result:
$$\frac{\#\K_{2i-1}(W_n(\bF_q);\bZ_p)}{\#\K_{2i-2}(W_n(\bF_q);\bZ_p)}=q^{i(n-1)}.$$
Both Brun and Angeltveit use classical trace methods and the $p$-adic filtration on the truncated Witt
vectors to translate part of the problem to the cases of truncated polynomial
rings where a complete answer is known.

The cases of $\K_3$ of $\bZ/p^n$ or $\bF_q[z]/z^2$ were also considered
earlier in~\cite{aisbett-lluis-puebla-snaith} using group homology calculations.

\section{New results}

As $\K(\Oscr/\varpi^k;\bZ_p)\we \tau_{\geq 0}\TC(\Oscr/\varpi^k;\bZ_p)$
by~\cite{dundas-goodwillie-mccarthy,hesselholt-madsen}, it is enough to
determine $\TC$ of these
rings. To do so, we use the filtration on $\TC$ constructed by
Bhatt--Morrow--Scholze in~\cite{bms2}. If $R$ is a quasisyntomic ring, there is
a complete decreasing filtration
$\F^{\geq \star}_\syn\TC(R;\bZ_p)$ with associated graded pieces
$$\F^{=i}_\syn\TC(R;\bZ_p)\we\bZ_p(i)(R)[2i],$$
where $\bZ_p(i)(R)$ is the weight $i$ syntomic cohomology of $R$
introduced in~\cite{bms2}. The syntomic complexes provide a $p$-adic analogue of the motivic
filtration on $\K$-theory.

As shown in~\cite{ammn}, the weight $i$ syntomic cohomology
$\bZ_p(i)(R)$ is concentrated in ${[0,i+1]}$, independent of $R$;
this means that $\H^r(\bZ_p(i)(R))=0$ for $r\notin[0,i+1]$. In the special case
of $\Oscr_K$ or $\Oscr_K/\varpi^n$, an argument using the $\varpi$-adic
associated graded implies that in fact the weight $i$ syntomic cohomology is in
${[0,2]}$; moreover, for $i\geq 1$,
$\H^0(\bZ_p(i)(\Oscr_K/\varpi^n))=0$ so the complex has cohomology concentrated
in degrees $1$ and $2$.

One checks that $\H^2(\bZ_p(1)(\Oscr_K/\varpi^n))=0$,
so the spectral sequence associated to the syntomic filtration
on $\TC$ collapses at the $\E_1$-page for $\Oscr/\varpi^n$ (or the $\E_2$-page
in the reindexing in~\cite[Thm.~1.12]{bms2}). Hence,
$$\TC_{2i-1}(\Oscr_K/\varpi^n;\bZ_p)\iso\H^1(\bZ_p(i)(\Oscr_K/\varpi^n))$$
for $i\geq 1$
and
$$\TC_{2i-2}(\Oscr_K/\varpi^n;\bZ_p)\iso\H^2(\bZ_p(i)(\Oscr_K/\varpi^n))$$
for $i\geq 2$. Thus, it makes sense to speak of the syntomic weights of
the $\K$-groups of $\Oscr_K/\varpi^n$.

\begin{theorem}\label{thm:main}
    For $i\geq 1$, if the residue field of $\Oscr_K$ has $q=p^f$ elements, then there is an explicit cochain complex
    \begin{gather*}
        \left(\bZ_p^{f(in-1)}\xrightarrow{\syn_0}\bZ_p^{2f(in-1)}\xrightarrow{\syn_1}\bZ_p^{f(in-1)}\right)
    \end{gather*}
    quasi-isomorphic to $\bZ_p(i)(\Oscr_K/\varpi^n)$.
    The terms are free $\bZ_p$-modules of the given ranks
     in cohomological degrees $0$, $1$, and $2$.
\end{theorem}

The proof of the existence of this explicit cochain complex model of the
syntomic complex will be discussed in Sections~\ref{sec:delta}
and~\ref{sec:syntomic}.

The groups $\K_*(\Oscr_K/\varpi^n)$ are torsion for $*>0$. In particular, the
complex above is exact rationally. Thus, to find the cohomology of
$\bZ_p(i)(\Oscr_K/\varpi^n)$, and hence the $p$-adic $\K$-groups of $\Oscr_K/\varpi^n$,
it is enough to compute the matrices $\syn_0$ and $\syn_1$ and their elementary
divisors.

\begin{theorem}
    The matrices $\syn_0$ and $\syn_1$ are effectively computable.
    Specifically, they can be determined with enough $p$-adic precision to
    guarantee computability of the effective divisors.
\end{theorem}

We have implemented our algorithm in {\ttfamily SAGE}~\cite{sagemath} in the case where $f=1$,
i.e., when the residue field is $\bF_p$. Future work will include an
implementation for general $f$.

\begin{corollary}
    There is an algorithm to determine the structure of $\K_r(\Oscr_K/\varpi^n)$ for any
    $K$, $n$, and $r$.
\end{corollary}

Along the way, we extend the result of Angeltveit on the quotients of the
orders from the unramified case to any $\Oscr_K/\varpi^n$.

\begin{corollary}\label{cor:quotient}
    For any $\Oscr_K/\varpi^n$,
    $$\frac{\#\K_{2i-1}(\Oscr_K/\varpi^n;\bZ_p)}{\#\K_{2i-2}(\Oscr_K/\varpi^n;\bZ_p)}=q^{i(n-1)},$$
    where $q=p^f$ is the order of the residue field of $\Oscr_K$.
\end{corollary}

This corollary is especially powerful thanks to the following theorem.

\begin{theorem}[Even vanishing theorem]\label{thm:even}
    If 
    \[
    i \geq \tfrac{p^2}{(p-1)^2} \big(p^{\lceil\tfrac{n}{e}\rceil}-1\big),
    \]
    then
    $\H^2(\bZ_p(i)(\Oscr_K/\varpi^n))=0$ and hence
    \[
    \K_{2i-2}(\Oscr_K/\varpi^n)=0
    \]
    if additionally $i\geq 2$.
\end{theorem}

\begin{corollary}\label{cor:odd}
    If   \[
    i \geq \tfrac{p^2}{(p-1)^2} \big(p^{\lceil\tfrac{n}{e}\rceil}-1\big),
    \]
    then $\#\K_{2i-1}(\Oscr_K/\varpi^k)=q^{i(n-1)}\cdot(q^i-1)$.
\end{corollary}

\begin{corollary}
    There is an algorithm to compute the orders of all of the $\K$-groups of
    $\Oscr/\varpi^n$.
\end{corollary}

Indeed, Theorem~\ref{thm:even} and Corollary~\ref{cor:odd} reduce the problem
to the computation of the cohomology of the syntomic complexes
$\bZ_p(i)(\Oscr/\varpi^n)$ for finitely many $i$: those satisfying
$$i<  \tfrac{p^2}{(p-1)^2} \left(p^{\lceil\tfrac{n}{e}\rceil}-1\right).
$$

This number grows rather quickly, but improvements are possible and will be described in our
forthcoming work~\cite{akn}.

\section{Computations}\label{sec:computations}

We present here four example calculations.

\subsection{$\bZ/4$}

The even vanishing theorem holds in syntomic weights $i\geq 12$. In fact, machine
computations show in this case that $\K_{2i-2}(\bZ/4)=0$ for all $i\geq 3$,
while $\K_2(\bZ/4)\iso\bZ/2$. Corollary~\ref{cor:quotient} together with
Quillen's calculation implies that
$$\#\K_3(\bZ/4)=8\cdot(2^2-1)\text{ and }\#\K_{2i-1}(\bZ/4)=2^i\cdot(2^i-1)$$ for $i\geq 3$.
This gives the complete calculation of the orders of all $\K$-groups of
$\bZ/4$.

The precise structure of the decomposition of $p$-primary part of the
$\K$-groups into cyclic groups remains unknown to us. Figure~\ref{fig:z4}
displays a table of the output of our machine computations giving the groups in syntomic weights
$i\leq 16$.

\begin{figure*}[h!]
    \centering
    \begin{tabular}{|r|r||r|r|} \hline
        $\K_1$ & $\bZ/2$ & $\K_{17}$ & $(\bZ/2)^3\oplus(\bZ/8)^2$\\
        $\K_2$ & $\bZ/2$&$\K_{18}$&$0$\\
        $\K_3$ & $\bZ/8$&$\K_{19}$&$\bZ/4\oplus\bZ/8\oplus\bZ/32$\\
        $\K_4$ & $0$&$\K_{20}$&$0$\\
        $\K_5$ & $\bZ/8$&$\K_{21}$&$(\bZ/2)^2\oplus(\bZ/4)^2\oplus\bZ/32$\\
        $\K_6$ & $0$&$\K_{22}$&$0$\\
        $\K_7$ & $\bZ/2\oplus\bZ/8$&$\K_{23}$&$(\bZ/2)^4\oplus\bZ/4\oplus\bZ/64$\\
        $\K_8$ & $0$&$\K_{24}$&$0$\\
        $\K_9$ &
        $(\bZ/2)^2\oplus\bZ/8$&$\K_{25}$&$(\bZ/2)^4\oplus\bZ/4\oplus\bZ/8\oplus\bZ/16$\\
        $\K_{10}$ & $0$&$\K_{26}$&$0$\\
        $\K_{11}$ &
        $\bZ/2\oplus\bZ/32$&$\K_{27}$&$\bZ/2\oplus\bZ/8\oplus\bZ/16\oplus\bZ/128$\\
        $\K_{12}$ & $0$&$\K_{28}$&$0$\\
        $\K_{13}$ &
        $\bZ/2\oplus\bZ/4\oplus\bZ/16$&$\K_{29}$&$(\bZ/2)^3\oplus(\bZ/4)^2\oplus\bZ/8\oplus\bZ/32$\\
        $\K_{14}$ & $0$&$\K_{30}$&$0$\\
        $\K_{15}$ &
        $(\bZ/2)^3\oplus\bZ/32$&$\K_{31}$&$(\bZ/2)^6\oplus\bZ/8\oplus\bZ/128$\\
        $\K_{16}$ & $0$&$\K_{32}$&$0$\\
        \hline
    \end{tabular}
    \caption{The $2$-adic $\K$-groups of $\bZ/4$ for syntomic weights $1\leq i\leq 16$;
    the final zero, $\K_{32}(\bZ/4;\bZ_2)=0$, is a (null) contribution from
    syntomic weight
    $17$.}
    \label{fig:z4}
\end{figure*}

\subsection{Chain rings of order $8$}

A chain ring is a commutative ring whose ideals are totally
ordered with respect to inclusion. Examples include valuation rings or quotients of valuation rings. Every
finite chain ring is of the form $\Oscr_K/\varpi^n$ for some $1\leq n<\infty$.
There are four chain rings of order $8$, namely $\bZ/8$,
$\bZ[2^{1/2}]/2^{3/2}$ (so $n=3$ in our notation), $\bF_2[z]/z^3$, and $\bF_8$;
see~\cite{clark-liang}. The $2$-adic $\K$-groups $\K_n(\bF_8;\bZ_2)$ vanish for
$n\geq 1$. Figure~\ref{fig:cr8} displays the low-degree $2$-adic $\K$-groups of the
other three chain rings of order $8$.

\begin{figure*}[h!]
    \centering
    \begin{tabular}{|r||r|r|r|}\hline
        $\K_r$&$\bZ/8$&$\bF_2[z]/z^3$&$\bZ_2[2^{1/2}]/2^{3/2}$\\
        \hline\hline
        $\K_1$&$\bZ/4$&$\bZ/4$&$\bZ/4$\\ \hline
        $\K_2$&$\bZ/2$&$0$&$0$\\\hline
        $\K_3$&$\bZ/4\oplus\bZ/8$&$\bZ/2\oplus\bZ/8$&$\bZ/2\oplus\bZ/8$\\\hline
        $\K_4$&$\bZ/2$&$0$&$0$\\\hline
        $\K_5$&$\bZ/2\oplus\bZ/64$&$(\bZ/2)^2\oplus\bZ/16$&$(\bZ/2)^2\oplus\bZ/16$\\\hline
        $\K_6$&$0$&$0$&$0$\\\hline
        $\K_7$&$(\bZ/4)^2$&$(\bZ/2)^2\oplus\bZ/4\oplus\bZ/16$&$(\bZ/2)^2\oplus\bZ/4\oplus\bZ/16$\\\hline
        $\K_8$&$0$&$0$&$0$\\\hline
        $\K_9$&$\bZ/2\oplus\bZ/4\oplus\bZ/128$&$(\bZ/2)^2\oplus(\bZ/4)^2\oplus\bZ/16$&$(\bZ/2)^2\oplus(\bZ/4)^2\oplus\bZ/16$\\\hline
        $\K_{10}$&$0$&$0$&$0$\\\hline
        $\K_{11}$&$\bZ/8\oplus\bZ/512$&$(\bZ/2)^3\oplus(\bZ/4)^2\oplus\bZ/32$&$(\bZ/2)^3\oplus(\bZ/4)^2\oplus\bZ/32$\\\hline
        $\K_{12}$&$0$&$0$&$0$\\\hline
        $\K_{13}$&$(\bZ/2)^2\oplus\bZ/8\oplus\bZ/512$&$(\bZ/2)^4\oplus\bZ/4\oplus\bZ/8\oplus\bZ/32$&$(\bZ/2)^4\oplus\bZ/4\oplus\bZ/8\oplus\bZ/32$\\\hline
        $\K_{14}$&$0$&$0$&$0$\\\hline
        $\K_{15}$&$(\bZ/2)^2\oplus\bZ/64\oplus\bZ/256$&$(\bZ/2)^4\oplus(\bZ/4)^2\oplus\bZ/8\oplus\bZ/32$&$(\bZ/2)^4\oplus(\bZ/4)^2\oplus\bZ/8\oplus\bZ/32$\\\hline
    \end{tabular}
    \caption{The $2$-adic $\K$-groups of the displayed chain rings of order $8$
    for syntomic weights $1\leq i\leq 8$. Note that the second and third
    columns agree.
    We do not know at present if this continues in all higher weights.
    The second column agrees with the calculations of~\cite{hesselholt-madsen-truncated} (see for
    example~\cite[Lem.~2]{speirs}).}
    \label{fig:cr8}
\end{figure*}

\subsection{Quotients of degree $2$ totally ramified $2$-adic fields}

The {\ttfamily lmfdb}~\cite{lmfdb} provides tables of $p$-adic fields based on
work of Jones--Roberts~\cite{jones-roberts}. There are $6$ totally ramified
degree $2$ extensions of $\bQ_2$. In Figure~\ref{fig:deg2}, we give
low-degree $p$-adic $\K$-groups of the quotients of these fields.

\begin{figure*}[h]
\centering
\begin{tikzpicture}
\node[scale=0.8] at (0,0) {\begin{tabular}{|r||r|r|r|r|r|r|r|} \hline
    {\ttfamily 2.2.2.1}&\multicolumn{7}{|l|}{$z^{2} + 2z + 2$}\\ \hline
$\K_{{r}}\backslash n$ & $\Oscr_K/\varpi^2$ & $\Oscr_K/\varpi^3$ & $\Oscr_K/\varpi^4$ & $\Oscr_K/\varpi^5$ & $\Oscr_K/\varpi^6$ & $\Oscr_K/\varpi^7$ & $\Oscr_K/\varpi^8$ \\ \hline \hline
$\K_{1}$ & $1$ & $2$ & $1, 2$ & $1, 1, 2$ & $1, 2, 2$ & $2, 2, 2$ & $2, 2, 3$ \\ \hline
$\K_{2}$ & $$ & $$ & $1$ & $1$ & $2$ & $2$ & $2$ \\ \hline
$\K_{3}$ & $1, 1$ & $1, 3$ & $1, 2, 4$ & $2, 3, 4$ & $1, 3, 3, 5$ & $1, 1, 3, 3, 6$ & $1, 1, 3, 4, 7$ \\ \hline
$\K_{4}$ & $$ & $$ & $$ & $$ & $1$ & $1$ & $2$ \\ \hline
$\K_{5}$ & $1, 1, 1$ & $1, 1, 4$ & $1, 2, 2, 4$ & $2, 2, 2, 6$ & $1, 2, 2, 4, 7$ & $1, 1, 2, 2, 4, 9$ & $1, 2, 2, 2, 6, 10$ \\ \hline
$\K_{6}$ & $$ & $$ & $$ & $$ & $$ & $1$ & $1$ \\ \hline
$\K_{7}$ & $1, 1, 1, 1$ & $1, 1, 2, 4$ & $1, 1, 1, 2, 3, 4$ & $1, 1, 1, 4, 4, 5$ & $1, 1, 1, 2, 5, 5, 5$ & $1, 1, 2, 3, 4, 7, 7$ & $1, 1, 1, 2, 3, 4, 8, 9$ \\ \hline
\end{tabular}};

\node[scale=0.8] at (0,-3.5) {\begin{tabular}{|r||r|r|r|r|r|r|r|} \hline
    {\ttfamily 2.2.2.2}&\multicolumn{7}{|l|}{$z^{2} + 2z - 2$}\\ \hline
$\K_{{r}}\backslash n$ & $\Oscr_K/\varpi^2$ & $\Oscr_K/\varpi^3$ & $\Oscr_K/\varpi^4$ & $\Oscr_K/\varpi^5$ & $\Oscr_K/\varpi^6$ & $\Oscr_K/\varpi^7$ & $\Oscr_K/\varpi^8$ \\ \hline \hline
$\K_{1}$ & $1$ & $2$ & $1, 2$ & $1, 1, 2$ & $1, 2, 2$ & $1, 2, 3$ & $1, 3, 3$ \\ \hline
$\K_{2}$ & $$ & $$ & $1$ & $1$ & $1$ & $1$ & $1$ \\ \hline
$\K_{3}$ & $1, 1$ & $1, 3$ & $1, 2, 4$ & $2, 3, 4$ & $1, 3, 3, 4$ & $1, 1, 3, 3, 5$ & $1, 1, 3, 4, 6$ \\ \hline
$\K_{4}$ & $$ & $$ & $$ & $$ & $1$ & $1$ & $2$ \\ \hline
$\K_{5}$ & $1, 1, 1$ & $1, 1, 4$ & $1, 2, 2, 4$ & $2, 2, 2, 6$ & $1, 2, 2, 4, 7$ & $1, 1, 2, 2, 4, 9$ & $1, 2, 2, 2, 6, 10$ \\ \hline
$\K_{6}$ & $$ & $$ & $$ & $$ & $$ & $1$ & $1$ \\ \hline
$\K_{7}$ & $1, 1, 1, 1$ & $1, 1, 2, 4$ & $1, 1, 1, 2, 3, 4$ & $1, 1, 1, 4, 4, 5$ & $1, 1, 1, 2, 5, 5, 5$ & $1, 1, 2, 3, 4, 7, 7$ & $1, 1, 1, 2, 3, 4, 8, 9$ \\ \hline
\end{tabular}};

\node[scale=0.8] at (0,-7) {\begin{tabular}{|r||r|r|r|r|r|r|r|}\hline
    {\ttfamily 2.2.3.1}& \multicolumn{7}{|l|}{$z^{2} + 14$}\\ \hline
$\K_{{r}}\backslash n$ & $\Oscr_K/\varpi^2$ & $\Oscr_K/\varpi^3$ & $\Oscr_K/\varpi^4$ & $\Oscr_K/\varpi^5$ & $\Oscr_K/\varpi^6$ & $\Oscr_K/\varpi^7$ & $\Oscr_K/\varpi^8$ \\ \hline \hline
$\K_{1}$ & $1$ & $2$ & $1, 2$ & $1, 1, 2$ & $1, 1, 3$ & $1, 2, 3$ & $1, 2, 4$ \\ \hline
$\K_{2}$ & $$ & $$ & $1$ & $1$ & $1$ & $1$ & $1$ \\ \hline
$\K_{3}$ & $1, 1$ & $1, 3$ & $1, 2, 4$ & $2, 3, 4$ & $1, 3, 3, 4$ & $1, 1, 3, 4, 4$ & $1, 2, 4, 4, 4$ \\ \hline
$\K_{4}$ & $$ & $$ & $$ & $$ & $1$ & $1$ & $2$ \\ \hline
$\K_{5}$ & $1, 1, 1$ & $1, 1, 4$ & $1, 2, 2, 4$ & $2, 2, 2, 6$ & $1, 2, 3, 3, 7$ & $1, 1, 2, 3, 3, 9$ & $1, 1, 3, 3, 5, 10$ \\ \hline
$\K_{6}$ & $$ & $$ & $$ & $$ & $$ & $1$ & $1$ \\ \hline
$\K_{7}$ & $1, 1, 1, 1$ & $1, 1, 2, 4$ & $1, 1, 2, 2, 2, 4$ & $1, 1, 2, 3, 4, 5$ & $1, 1, 2, 2, 4, 5, 5$ & $1, 2, 2, 3, 5, 5, 7$ & $1, 1, 2, 2, 3, 5, 6, 9$ \\ \hline
\end{tabular}};

\node[scale=0.8] at (0,-10.5) {\begin{tabular}{|r||r|r|r|r|r|r|r|} \hline
    {\ttfamily 2.2.3.2}& \multicolumn{7}{|l|}{$z^{2} + 6$}\\ \hline
$\K_{{r}}\backslash n$ & $\Oscr_K/\varpi^2$ & $\Oscr_K/\varpi^3$ & $\Oscr_K/\varpi^4$ & $\Oscr_K/\varpi^5$ & $\Oscr_K/\varpi^6$ & $\Oscr_K/\varpi^7$ & $\Oscr_K/\varpi^8$ \\ \hline \hline
$\K_{1}$ & $1$ & $2$ & $1, 2$ & $1, 1, 2$ & $1, 1, 3$ & $1, 2, 3$ & $1, 2, 4$ \\ \hline
$\K_{2}$ & $$ & $$ & $1$ & $1$ & $1$ & $1$ & $1$ \\ \hline
$\K_{3}$ & $1, 1$ & $1, 3$ & $1, 2, 4$ & $2, 3, 4$ & $1, 3, 3, 4$ & $1, 1, 3, 4, 4$ & $1, 2, 3, 4, 5$ \\ \hline
$\K_{4}$ & $$ & $$ & $$ & $$ & $1$ & $1$ & $2$ \\ \hline
$\K_{5}$ & $1, 1, 1$ & $1, 1, 4$ & $1, 2, 2, 4$ & $2, 2, 2, 6$ & $1, 2, 3, 3, 7$ & $1, 1, 2, 3, 3, 9$ & $1, 1, 3, 3, 5, 10$ \\ \hline
$\K_{6}$ & $$ & $$ & $$ & $$ & $$ & $1$ & $1$ \\ \hline
$\K_{7}$ & $1, 1, 1, 1$ & $1, 1, 2, 4$ & $1, 1, 2, 2, 2, 4$ & $1, 1, 2, 3, 4, 5$ & $1, 1, 2, 2, 4, 5, 5$ & $1, 2, 2, 3, 5, 5, 7$ & $1, 1, 2, 2, 3, 5, 6, 9$ \\ \hline
\end{tabular}};

\node[scale=0.8] at (0,-14) {\begin{tabular}{|r||r|r|r|r|r|r|r|} \hline
    {\ttfamily 2.2.3.3}& \multicolumn{7}{|l|}{$z^{2} + 2$}\\ \hline
$\K_{{r}}\backslash n$ & $\Oscr_K/\varpi^2$ & $\Oscr_K/\varpi^3$ & $\Oscr_K/\varpi^4$ & $\Oscr_K/\varpi^5$ & $\Oscr_K/\varpi^6$ & $\Oscr_K/\varpi^7$ & $\Oscr_K/\varpi^8$ \\ \hline \hline
$\K_{1}$ & $1$ & $2$ & $1, 2$ & $1, 1, 2$ & $1, 1, 3$ & $1, 2, 3$ & $1, 2, 4$ \\ \hline
$\K_{2}$ & $$ & $$ & $1$ & $1$ & $1$ & $1$ & $1$ \\ \hline
$\K_{3}$ & $1, 1$ & $1, 3$ & $1, 2, 4$ & $2, 3, 4$ & $1, 3, 3, 4$ & $1, 1, 3, 4, 4$ & $1, 2, 3, 4, 5$ \\ \hline
$\K_{4}$ & $$ & $$ & $$ & $$ & $1$ & $1$ & $2$ \\ \hline
$\K_{5}$ & $1, 1, 1$ & $1, 1, 4$ & $1, 2, 2, 4$ & $2, 2, 2, 6$ & $1, 2, 3, 3, 7$ & $1, 1, 2, 3, 3, 9$ & $1, 1, 3, 3, 5, 10$ \\ \hline
$\K_{6}$ & $$ & $$ & $$ & $$ & $$ & $1$ & $1$ \\ \hline
$\K_{7}$ & $1, 1, 1, 1$ & $1, 1, 2, 4$ & $1, 1, 2, 2, 2, 4$ & $1, 1, 2, 3, 4, 5$ & $1, 1, 2, 2, 4, 5, 5$ & $1, 2, 2, 3, 5, 5, 7$ & $1, 1, 2, 2, 3, 5, 6, 9$ \\ \hline
\end{tabular}};

\node[scale=0.8] at (0,-17.5) {\begin{tabular}{|r||r|r|r|r|r|r|r|} \hline
    {\ttfamily 2.2.3.4}& \multicolumn{7}{|l|}{$z^{2} + 10$}\\ \hline
$\K_{{r}}\backslash n$ & $\Oscr_K/\varpi^2$ & $\Oscr_K/\varpi^3$ & $\Oscr_K/\varpi^4$ & $\Oscr_K/\varpi^5$ & $\Oscr_K/\varpi^6$ & $\Oscr_K/\varpi^7$ & $\Oscr_K/\varpi^8$ \\ \hline \hline
$\K_{1}$ & $1$ & $2$ & $1, 2$ & $1, 1, 2$ & $1, 1, 3$ & $1, 2, 3$ & $1, 2, 4$ \\ \hline
$\K_{2}$ & $$ & $$ & $1$ & $1$ & $1$ & $1$ & $1$ \\ \hline
$\K_{3}$ & $1, 1$ & $1, 3$ & $1, 2, 4$ & $2, 3, 4$ & $1, 3, 3, 4$ & $1, 1, 3, 4, 4$ & $1, 2, 3, 4, 5$ \\ \hline
$\K_{4}$ & $$ & $$ & $$ & $$ & $1$ & $1$ & $2$ \\ \hline
$\K_{5}$ & $1, 1, 1$ & $1, 1, 4$ & $1, 2, 2, 4$ & $2, 2, 2, 6$ & $1, 2, 3, 3, 7$ & $1, 1, 2, 3, 3, 9$ & $1, 1, 3, 3, 5, 10$ \\ \hline
$\K_{6}$ & $$ & $$ & $$ & $$ & $$ & $1$ & $1$ \\ \hline
$\K_{7}$ & $1, 1, 1, 1$ & $1, 1, 2, 4$ & $1, 1, 2, 2, 2, 4$ & $1, 1, 2, 3, 4, 5$ & $1, 1, 2, 2, 4, 5, 5$ & $1, 2, 2, 3, 5, 5, 7$ & $1, 1, 2, 2, 3, 5, 6, 9$ \\ \hline
\end{tabular}};
\end{tikzpicture}

\caption{The $2$-adic $\K$-groups in syntomic weights $i=1,2,3,4$ for the totally
ramified degree $2$ extensions of $\bZ_2$. The {\ttfamily lmfdb}~\cite{lmfdb} label is given
in the top left corner together with an Eisenstein polynomial. The data gives the
exponents of the elementary divisors in each degree: for example, the entry
$1,3$ in the $\K_3$ row of the $\Oscr_K/\varpi^3$ column means that
$\K_3(\Oscr_K/\varpi^3;\bZ_2)\iso\bZ/2\oplus\bZ/8$.}
\label{fig:deg2}
\end{figure*}

\subsection{$\bZ/9$}

The even vanishing theorem holds in syntomic weights $i\geq 18$.
Figure~\ref{fig:z9} displays a table of
the output of our machine computations in syntomic weights $i\leq 18$.
In particular, $\K_4(\bZ/9)\iso\bZ/3$ and all other positive even
$\K$-groups vanish. In odd degrees,
$$\#\K_5(\bZ/9)=81\cdot(3^3-1)\text{ and }\#\K_{2i-1}(\bZ/9)=3^i\cdot(3^i-1)$$
for $i\geq 1$, $i\neq 3$.
This gives the complete calculation of the orders of all $\K$-groups of
$\bZ/9$.

\begin{figure*}[h!]
    \centering
    \begin{tabular}{|r|r||r|r|} \hline
        $\K_1$ & $\bZ/3$ &$\K_{19}$&$(\bZ/3)^3\oplus\bZ/9\oplus\bZ/243$\\
        $\K_2$ & $0$& $\K_{20}$ & $0$\\
        $\K_3$ & $(\bZ/3)^2$& $\K_{21}$ & $(\bZ/3)^3\oplus\bZ/9\oplus\bZ/729$\\
        $\K_4$ & $\bZ/3$& $\K_{22}$ & $0$\\
        $\K_5$ & $\bZ/81$& $\K_{23}$ & $\bZ/3\oplus\bZ/27\oplus\bZ/6561$\\
        $\K_6$ & $0$& $\K_{24}$ &$0$\\
        $\K_7$ & $\bZ/3\oplus\bZ/27$& $\K_{25}$ & $(\bZ/3)^4\oplus\bZ/9\oplus\bZ/2187$\\
        $\K_8$ & $0$ & $\K_{26}$ & $0$\\
        $\K_9$ & $\bZ/3\oplus\bZ/81$ & $\K_{27}$ & $(\bZ/3)^4\oplus\bZ/9\oplus\bZ/6561$\\
        $\K_{10}$ & $0$ & $\K_{28}$ & $0$\\
        $\K_{11}$ & $(\bZ/27)^2$ & $\K_{29}$ & $\bZ/3\oplus\bZ/9\oplus\bZ/27\oplus\bZ/19683$\\
        $\K_{12}$ & $0$ & $\K_{30}$ & $0$\\
        $\K_{13}$ & $(\bZ/3)^2\oplus\bZ/243$ & $\K_{31}$ & $(\bZ/3)^4\oplus(\bZ/9)^2\oplus\bZ/6561$\\
        $\K_{14}$ & $0$ & $\K_{32}$ & $0$\\
        $\K_{15}$& $(\bZ/3)^2\oplus\bZ/729$ & $\K_{33}$ & $(\bZ/3)^4\oplus(\bZ/9)^2\oplus\bZ/19683$\\
        $\K_{16}$&$0$ & $\K_{34}$ & $0$\\
        $\K_{17}$&$\bZ/9\oplus\bZ/2187$ & $\K_{35}$ & $(\bZ/9)^2\oplus\bZ/243\oplus\bZ/19683$\\
        $\K_{18}$&$0$ & $\K_{36}$ & $0$\\

        \hline
    \end{tabular}
    \caption{The $3$-adic $\K$-groups of $\bZ/9$ for syntomic weights $1\leq i\leq 18$.
    The contribution of $\K_{36}(\bZ/9;\bZ_3)=0$ is a (null) group from weight $19$.}
    \label{fig:z9}
\end{figure*}

\section{Prismatic cohomology over $\delta$-rings}\label{sec:delta}

Our proofs are motivated by previous work of
Krause--Nikolaus~\cite{krause-nikolaus} and the approach of
Liu--Wang~\cite{liu-wang}. There are two main new ideas: the notion
of prismatic cohomology relative to a $\delta$-ring and
the systematic use of the filtration on the syntomic complexes induced by the $\varpi$-adic
filtration on $\Oscr_K/\varpi^n$. Similar filtrations have also been used  by Angeltveit~\cite{angeltveit} and Brun~\cite{brun} in the topological context.

Let $A^0=W(\bF_q)\llbracket z\rrbracket$ be the $\delta$-ring with $\delta(z)=0$ and
hence $\varphi(z)=z^p$. If $E(z)$ is an Eisenstein polynomial for $\Oscr_K$, then
the pair $(A^0,(E(z)))$ is a prism. Bhatt and Scholze show that $\Prism_{(\Oscr_K/\varpi^n)/A^0}$ is discrete
and admits a description as a prismatic envelope $A^0\{\frac{\varpi^n}{E(z)}\}^\wedge$
in the sense of~\cite[Prop.~3.13]{prisms};
the prismatic envelope is an explicit pushout in
$(p,E(z))$-complete $\delta$-rings over $A^0$.

The main idea is to determine the syntomic complexes $\bZ_p(i)(\Oscr/\varpi^n)$
by descent along the map
$\Prism_{\Oscr/\varpi^n}\rightarrow\Prism_{(\Oscr/\varpi^n)/A^0}$ from absolute
prismatic cohomology to relative prismatic cohomology.
To make sense of this, we introduce prismatic cohomology relative to a
$\delta$-ring. Let us outline the definition.

Given an arbitrary derived $p$-complete $\delta$-ring $A$ and a derived
$p$-complete $A$-algebra $R$, let $X=\Spf R$ and let $(X/A)_\Prism$ be the
opposite of the category of commutative diagrams
$$\xymatrix{
    A\ar[r]\ar[d]&B\ar[d]\\
    R\ar[r]&B/J,
}$$ where $(B,J)$ is a bounded prism and $A\rightarrow B$ is a map of
$\delta$-rings.

By definition, $\Prism_{R/A}=\R\Gamma((X/A)_\Prism,\Oscr_\Prism)$, where
$\Oscr_\Prism$ is
the prismatic structure sheaf, which sends a commutative diagram as above to $B$.
Warning: this site-theoretic definition should be derived in general, but gives
the correct answer under additional assumptions on $R$, in particular in the case of $R=\Oscr_K/\varpi^n$ over the multivariable
Breuil--Kisin prisms appearing in this paper. 

\begin{example}
    If $A=\bZ_p$ is the initial (derived $p$-complete) $\delta$-ring, then
    $\Prism_{R/\bZ_p}$ recovers absolute prismatic cohomology as introduced
    in~\cite{bms2,prisms} and studied further in~\cite{bhatt-lurie-apc}.
    More generally, this is true if $A$ is replaced by the ring of $p$-typical
    Witt vectors of any perfect $\bF_p$-algebra.
\end{example}

\begin{example}\label{ex:prismatic}
    If $(A,I)$ is a prism and $R$ is an $A/I$-algebra, then $\Prism_{R/A}$
    agrees with derived relative prismatic cohomology as studied in~\cite{prisms}.
\end{example}


Now, consider the augmented cosimplicial diagram $A^\bullet$ where
$A^{-1}=W(\bF_q)$, $A^0=W(\bF_q)\llbracket z\rrbracket$, and
$A^s=W(\bF_q)\llbracket z_0,\ldots,z_s\rrbracket$. This is a completed descent
complex for $W(\bF_q)\rightarrow W(\bF_q)[z]$.

In the cosimplicial diagram
$$W(\bF_q)\stack{1} A^0\stack{3} A^1\stack{5} A^2\cdots,$$
the arrows are all $\delta$-ring maps and the entire diagram admits a
map to $\Oscr_K$ sending each generator $z_j$ to $\varpi$.
As a result, for any $\Oscr_K$-algebra $R$, there is an induced augmented cosimplicial
diagram
in prismatic cohomology of $R$ relative to the $\delta$-rings $A^\bullet$.

\begin{theorem}\label{thm:cosimplicial}
    The augmented cosimplicial diagram
$$\Prism_{R}\stack{1} \Prism_{R/A^0}\stack{3}
\Prism_{R/A^1}\stack{5}\Prism_{R/A^2}\cdots$$
is a limit diagram for $R=\Oscr_K/\varpi^n$.
\end{theorem}

Thus, the absolute prismatic cohomology of an $\Oscr_K$-algebra, such as
$\Oscr_K/\varpi^n$, can be computed by descent using the cosimplicial diagram
above.

This does not make sense when speaking of prismatic cohomology as defined
in~\cite{prisms} because there is no compatible way to equip the entire
cosimplicial diagram with the structure of a cosimplicial prism. For example,
if $E(z)$ is an Eisenstein polynomial making $A^0=W(\bF_q)\llbracket
z\rrbracket$ into a prism, both $E(z_0)$
and $E(z_1)$ are distinguished elements in $A^1=W(\bF_q)\llbracket
z_0,z_1\rrbracket$ making it into a prism in two different ways.

\begin{proposition}\label{prop:envelope}
    For any $s\geq 0$, the relative prismatic cohomology
    $\Prism_{(\Oscr_K/\varpi^n)/A^s}$ is discrete and is isomorphic to a prismatic envelope
    $$A^s\left\{
        \frac{z_0^n}{E(z_0)},\frac{z_1-z_0}{E(z_0)},\ldots,\frac{z_n-z_0}{E(z_0)}\right\}^\wedge.$$
\end{proposition}

The proposition follows immediately from Example~\ref{ex:prismatic}.
Note that while prismatic cohomology relative to $\delta$-rings is functorial in
arbitrary maps of $\delta$-rings, the presentation of a given term
$\Prism_{R/A^s}$ as a prismatic envelope depends on the choice of a prism
structure $J$ on $A^s$ making $R$ into an $A^s/J$-algebra. In the theorem
above, we choose to make $A^s$ into a prism with respect to the ideal
$(E(z_0))$.

It follows that the cosimplicial diagram appearing in
Theorem~\ref{thm:cosimplicial} gives a resolution of
$\Prism_{\Oscr_K/\varpi^n}$ as the limit of a cosimplicial diagram of discrete
$\delta$-rings.

To give the main idea of the rest of the argument, we illustrate it here for
prismatic cohomology instead of the syntomic complexes. The absolute prismatic cohomology
of a quasisyntomic ring $R$
admits a Nygaard filtration $\Nscr^{\geq\star}\Prism_R$; Nygaard completion of
prismatic cohomology is written $\Prismhat_R$.

\begin{proposition}
    The Nygaard-complete absolute prismatic cohomology groups $\H^r(\Prismhat_{\Oscr_K/\varpi^n})$
    vanish for $r\neq 0,1$.
\end{proposition}

The proposition can be proved by computing directly with a Nygaard-complete,
Frobenius-twisted
variant of the cosimplicial diagram in Theorem~\ref{thm:cosimplicial} using the prismatic
envelopes of
Proposition~\ref{prop:envelope}. Alternatively, one can argue 
as follows: the $\varpi$-adic filtration on $\Oscr_K/\varpi^n$ induces a
filtration on $\Prism_{\Oscr_K/\varpi^n}$ whose completion agrees with $\Prismhat_{\Oscr_K/\varpi^n}$, and whose associated graded is the same as
that of the corresponding filtration on $\Prismhat_{\bF_q[z]/z^n}$. This
associated graded can be described using crystalline cohomology and
vanishes away from cohomological degrees $0,1$. Thus, by d\'evissage and
completeness, the same vanishing holds for $\Prismhat_{\Oscr_K/\varpi^n}$.

It follows from the proposition that the cochain complex $A^0\rightarrow
A^1\rightarrow A^2\rightarrow\cdots$  associated to
the cosimplicial abelian group $\Prismhat_{(\Oscr_K/\varpi^n)/A^\bullet}$ is exact in degrees $\geq 2$.
This reduces the computation of $\Prismhat_{\Oscr_K/\varpi^n}$ to a much smaller
computation involving prismatic envelopes of $\Oscr_K/\varpi^n$ relative to
$A^0$, $A^1$, and $A^2$.

However, we are interested not in the absolute prismatic cohomology of
$\Oscr_K/\varpi^n$ but rather in its syntomic cohomology. Relative syntomic
cohomology is defined in the setting of prismatic cohomology relative to a
$\delta$-ring. We first have to explain the Nygaard filtration
and the Breuil--Kisin twist, following~\cite{prisms,bhatt-lurie-apc}.

The Frobenius twist $\Prism_{R/A}^{(1)}$ is defined
to be $\Prism_{R/A}\otimes_A{_\varphi A}$, the base-change of $\Prism_{R/A}$
along the Frobenius map on $A$. The Frobenius twist admits a map
$\Prism_{R/A}^{(1)}\rightarrow\Prism_{R/A}$ and the Nygaard filtration
$\Nscr^{\geq\star}\Prism_{R/A}^{(1)}$ is
a filtration which is taken by this map to the $I$-adic filtration on
$\Prism_{R/A}$. If $\Prism_{R/A}$ is discrete (as in our examples of interest)
then the Nygaard filtration is simply the preimage of the $I$-adic filtration.

Given a prism $(A,I)$, let $I_r$ be the invertible $A$-module
$I\cdot\varphi(I)\cdots\varphi^{r-1}(I)$. If $(A,I)$ is transversal, meaning
that $A/I$ is $p$-torsion-free, then the
canonical map $I_r/I_r^2\rightarrow I_{r-1}/I_{r-1}^2$ is divisible by $p$ and
the induced map $I_r/I_r^2\xrightarrow{1/p}I_{r-1}/I_{r-1}^2$ is surjective.
The Breuil--Kisin twist is defined to be
$$A\{1\}=\lim\left(\cdots\rightarrow
I_3/I_3^2\xrightarrow{1/p}I_2/I_2^2\xrightarrow{1/p}I/I^2\right).$$
This is an invertible $A$-module.
For a general $A$-module $M$, let $M\{1\}=M\otimes_A A\{1\}$.

The relative syntomic cohomology of $R$ over a $\delta$-ring $A$ is
$$\bZ_p(i)(R/A)=\mathrm{fib}\left(\Nscr^{\geq
i}\Prism^{(1)}_{R/A}\{i\}\xrightarrow{\can-\varphi}\Prism^{(1)}_{R/A}\{i\}\right),$$
where $\varphi$ is a Frobenius which exists on $\Nscr^{\geq
i}\Prism_{R/A}^{(1)}\{i\}$.
Note that in~\cite{bms2}, the syntomic complexes are defined using Nygaard
complete prismatic cohomology; however, the two definitions agree
by~\cite[Lem.~7.22]{bms2} or~\cite[Cor.~5.31]{ammn}.

It follows along the lines of Theorem~\ref{thm:cosimplicial} that, for each
$i\geq 0$, the limit of the cosimplicial diagram
$$\bZ_p(i)(R/A^0)\stack{3}\bZ_p(i)(R/A^1)\stack{5}\cdots$$
is equivalent to $\bZ_p(i)(R)$ when $R=\Oscr_K/\varpi^k$.

The fact that the Nygaard-complete absolute prismatic cohomology $\Prismhat_{\Oscr_K/\varpi^n}$ is
concentrated in cohomological degrees $0,1$ implies that
$\bZ_p(i)(\Oscr_K/\varpi^n)$ is concentrated in cohomological degrees $0,1,2$.
In fact, it is not hard to show that, for $i\geq 1$, each relative syntomic
complex $\bZ_p(i)((\Oscr_K/\varpi^n)/A^s)$ is concentrated in cohomological degree $1$.
Thus, the spectral sequence associated to the limit diagram
$$\bZ_p(i)(\Oscr_K/\varpi^n)\we\lim_\Delta\bZ_p(i)((\Oscr_K/\varpi^n)/A^\bullet)$$
implies that $\bZ_p(i)(\Oscr_K/\varpi^n)$ is concentrated in cohomological
degrees $1,2$ for $i\geq 1$.

By the same spectral sequence, to determine $\bZ_p(i)(\Oscr_K/\varpi^n)$,
and hence $\K_{2i-2}(\Oscr_K/\varpi^n;\bZ_p)$ and
$\K_{2i-1}(\Oscr_K/\varpi^n;\bZ_p)$, it is enough to compute the cohomology of
the complex
\begin{gather*}
    \H^1(\bZ_p(i)(R/A^0)\\\rightarrow\ker\left(\H^1(\bZ_p(i)(R/A^1))\rightarrow\H^1(\bZ_p(i)(R/A^2))\right)
\end{gather*}
where $R=\Oscr_K/\varpi^n$.
In the next section, we explain how to use the $\varpi$-adic filtration to
reduce this to a finite problem.

\section{The syntomic matrices}\label{sec:syntomic}

In the cosimplicial diagram $A^\bullet$, each term is a filtered $\delta$-ring,
where in $A^s=W(k)\llbracket z_0,\ldots,z_s\rrbracket$ the weight of $z_j$ is
$1$. A filtered $\delta$-ring is a $\delta$-ring $A$ with a complete and
separated decreasing
filtration $\Fscr^{\geq \star} A$ such that $\delta(\Fscr^{\geq i}A)\subseteq\Fscr^{\geq pi}A$. Since each
$A^\bullet\rightarrow\Oscr_K/\varpi^n$ is a filtered map where
$\Oscr_K/\varpi^n$ is given the $\varpi$-adic filtration, all resulting
invariants, such as prismatic or syntomic cohomology complexes admit induced
filtrations, which we will write for instance as
$\Fscr^{\geq \star}\bZ_p(i)((\Oscr_K/\varpi^n)/A^\bullet)$.

\begin{theorem}
    For $b\geq in-1$ and $i\geq 1$, the natural maps
    $$\xymatrix{
        &\Fscr^{[1,b]}\bZ_p(i)(\Oscr_K/\varpi^n)\ar[d]\\
        \bZ_p(i)(\Oscr_K/\varpi^n)\ar[r]&\Fscr^{[0,b]}\bZ_p(i)(\Oscr_K/\varpi^n)
    }$$
    are equivalences.
\end{theorem}

The right-hand arrow is easy to handle because
$\Fscr^{=0}\bZ_p(i)(\Oscr_K/\varpi^n)\we\bZ_p(i)(\bF_q)\we 0$ for $i>0$. 
For the left-hand arrow, we argue by an explicit study of the interaction
between the $\Fscr$-filtration and the Nygaard filtration on each
$\Prism_{(\Oscr_K/\varpi^n)/A^\bullet}$.

The entire problem has now been reduced to a finite computation.
Set $R=\Oscr_K/\varpi^n$ and consider the commutative diagram
$$\xymatrix{
    \Fscr^{[1,b]}\Nscr^{\geq
    i}\Prism^{(1)}_{R/A^0}\{i\}\ar[r]\ar[d]&\Fscr^{[1,b]}\Nscr^{\geq i}\Prism^{(1)}_{R/A^1}\{i\}\ar[d]\\
    \Fscr^{[1,b]}\Prism^{(1)}_{R/A^0}\{i\}\ar[r]&\Fscr^{[1,b]}\Prism^{(1)}_{R/A^1}\{i\}.
}$$
All four terms are finitely generated free $\bZ_p$-modules.
The vertical fibers are $\bZ_p(i)(R/A^0)$ and $\bZ_p(i)(R/A^1)$, respectively.

Our approach to the computation avoids the more traditional approach
of computing either $\TR(\Oscr_K/\varpi^n)^{F=1}$ or computing
$\TC(\Oscr_K/\varpi^n)$ as the fiber of
$\TC^-(\Oscr_K/\varpi^n)\xrightarrow{\can-\varpi}\TP(\Oscr_K/\varpi^n)$. It
would nevertheless be very interesting to understand $\TP(\Oscr_K/\varpi^n)$.

Since the complexes $\Fscr^{[1,b]}\Nscr^{\geq i}\Prism_{R}\{i\}$ and
$\Fscr^{[1,b]}\Prism_R\{i\}$ are torsion
for $i\geq 1$ by another use of the $\varpi$-adic filtration, one can replace
$$\ker\left(\Fscr^{[1,b]}\Nscr^{\geq
i}\Prism^{(1)}_{R/A^1}\{i\}\rightarrow\Fscr^{[1,b]}\Nscr^{\geq
i}\Prism^{(1)}_{R/A^2}\{i\}\right)$$ with the saturation of the image of the top
horizontal map, where by saturation we mean the sub-$\bZ_p$-module consisting
of elements $x$ such that $p^Nx$ is in the image for some $N$, and similarly
for
$\ker\left(\Fscr^{[1,b]}\Prism^{(1)}_{R/A^1}\{i\}\rightarrow\Fscr^{[1,b]}\Prism^{(1)}_{R/A^2}\{i\}\right)$.
Write $S^0$ and $S^1$ for the saturations.
The resulting commutative square
$$\xymatrix{
    \Fscr^{[1,b]}\Nscr^{\geq i}\Prism^{(1)}_{R/A^0}\{i\}\ar[r]\ar[d]&S^0\ar[d]\\
    \Fscr^{[1,b]}\Prism^{(1)}_{R/A^0}\{i\}\ar[r]&S^1
}$$
consists of free $\bZ_p$-modules of rank $bf$ and the total cohomology computes
$\Fscr^{[1,b]}\bZ_p(i)(R)$ and hence $\bZ_p(i)(R)=\bZ_p(i)(\Oscr_K/\varpi^n)$ for
$i\geq 1$.

To conclude, we use explicit polynomial presentations of the relevant prismatic
envelopes as well as Breuil--Kisin
orientations to give explicit bases of all four terms and to compute the maps between them.
Taking $b=in-1$, the result is the matrices $\syn_0$ and $\syn_1$ and the complex appearing in
Theorem~\ref{thm:main}.

\paragraph{Acknowledgements.}
We are very grateful to Bhargav Bhatt, Lars Hesselholt, Akhil Mathew, Noah Riggenbach, Peter Scholze, and
Chuck Weibel for comments on a draft of this announcement.

The first author was supported by NSF grants
DMS-2102010 and DMS-2120005 and a Simons Fellowship; he would like to thank
Universit\"at M\"unster
for its hospitality during a visit in 2020. 
The second and third author were funded by the Deutsche Forschungsgemeinschaft
(DFG, German Research Foundation) – Project-ID 427320536 – SFB 1442, as well as
under Germany’s Excellence Strategy EXC 2044 390685587, Mathematics Münster:
Dynamics–Geometry–Structure. They would also like to thank the Mittag--Leffler
Institute for its hospitality while working on this project.

\small
\bibliographystyle{amsplain}
\bibliography{kzpn}

\printindex

\end{document}